%% file: inducedsubgraph.tex
\documentclass[11pt,reqno]{amsart}
\usepackage[top=1in,bottom=1in]{geometry}
\usepackage{amsmath}
\usepackage{amssymb}
\usepackage{amsfonts}
\usepackage{amscd}
\usepackage{amsthm}
\usepackage{amsxtra}
\usepackage{stmaryrd}

\usepackage[all]{xy}  \CompileMatrices
\usepackage{mathrsfs}
\usepackage{nicefrac}
\usepackage{graphicx}
\usepackage{color}
\usepackage[raggedright,IT,hang]{subfigure}
\usepackage{algorithmic}

\input{defs.amsart}

\DeclareMathOperator{\cop}{cop}
\DeclareMathOperator{\tw}{tw}
\DeclareMathOperator{\circum}{circum}

\newcounter{myremembertheorem}

\newcommand{\defn}[1]{\textit{#1}}

\begin{document}

\title{The Cops and Robber game on graphs with forbidden (induced)~subgraphs}%
\thanks{This work was partially supported by the
\textit{Actions de Recherche Concert\'ees (ARC)\,}
fund of the \textit{Communaut\'e fran\c{c}aise de Belgique}.}
\author{Gwena\"el Joret}
\address{Gwena\"el Joret\\
  D\'epartement d'Informatique\\
  Universit\'e Libre de Bruxelles\\
  Brussels, Belgium}%
\email{gjoret@ulb.ac.be}
\thanks{G.~Joret is a Research Fellow of the \textit{Fonds
National de la Recherche Scientifique (F.R.S.--FNRS)}}
\author{Marcin Kami\'nski}
\address{Marcin Kami\'nski\\
  D\'epartement d'Informatique\\
  Universit\'e Libre de Bruxelles\\
  Brussels, Belgium}%
\email{Marcin.Kaminski@ulb.ac.be}
\author{Dirk Oliver Theis}
\address{Dirk Oliver Theis\\
  Service de G\'eometrie Combinatoire et Th\'eorie des Groupes\\
  D\'epartement de Math\'ematique\\
  Universit\'e Libre de Bruxelles\\
  Brussels, Belgium}%
\email{Dirk.Theis@ulb.ac.be}%
\thanks{M.~Kami\'nski and 
D.O.~Theis supported by \textit{Fonds National de la Recherche Scientifique
(F.R.S.--FNRS)}}


\date{\today}


\begin{abstract}
  The two-player, complete information game of Cops and Robber is played on undirected finite graphs. A number of cops and one robber are positioned on vertices
  and take turns in sliding along edges. The cops win if, after a move, a cop and the robber are on the same vertex.  The minimum number of cops needed to catch
  the robber on a graph is called the cop number of that graph.

  In this paper, we study the cop number in the classes of graphs defined by forbidding one or more graphs as either subgraphs or induced subgraphs. In the case of
  a single forbidden graph we completely characterize (for both relations) the graphs which force bounded cop number. 
  En passant, we bound the cop number in terms of the tree-width.
\end{abstract}
\maketitle

\section{Introduction}

Graphs studied in this paper are finite, undirected, without loops and multiple edges. We use standard notation and terminology; for what is not defined here, we
refer the reader to Diestel~\cite{Diestel-GTIII}.

The game of \defn{Cops and Robber} is played on a connected graph by two players -- the cops and the robber. The cop player has at her disposal $k$ pieces (cops),
for some integer $k \geq 1$, and the cop player has only one piece (the robber).  
The game begins with
the cop player placing her $k$ cops on (not necessarily distinct) vertices of the graph.  Next, the robber player chooses a vertex for his piece. Now, starting
with the cop player, the two players move their pieces alternately. In the cops' move, she decides for each of her cops whether it stands still or is moved to an
adjacent vertex. In the robber's move, he can choose to move or not to move the piece. The game ends when a cop and the robber are on the same vertex (that is, the
cops catch the robber); in this case the cop player wins. The robber wins if he can never be caught by the cops. Both players have complete information, that is,
they know the graph and the positions of all the pieces.

The key problem is to know how many cops are needed to catch a robber on a given graph. For a connected graph $G$, the smallest integer $k$ such that with $k$
cops, the cop player has a winning strategy is called the \defn{cop number} of $G$ and is denoted by $\cop(G)$. We follow Berarducci and
Intrigila~\cite{BerarducciIntrigila93} in defining the cop number of a non-connected graph as the maximum cop number of its connected components.  Nowakowski and
Winkler \cite{NowakowskiWinkler83} and Quilliot \cite{QuilliotTh} have characterized the class of graphs with cop number 1.  Finding a~combinatorial
characterization of graphs with cop number $k$ (for $k\ge 2$) is a major open problem in the field (see \cite{ClarkePhD}).  On the other hand, algorithmic
characterizations of such graphs, which are polynomial in the size of the graph but not in $k$, do exist~\cite{BerarducciIntrigila93,GoldsteinReingold95}.
However, determining the cop number of a graph is a computationally hard problem \cite{FominComplexity}.  The cop number of random graphs has been studied
\cite{BollobasKunLeader08,BonaHahnWang07}.  For literature review we refer the reader to a recent survey on graph searching \cite{SurveyFomThi}.

In this paper we study the cop number for different types of graph classes. Our motivation is to learn what structural properties of graphs force the cop number to
be bounded. (We say that the cop number is \defn{bounded} for a class of graphs, if there exists a constant $C$ such that the cop number of every graph from the
class is at most $C$; otherwise the cop number is \defn{unbounded} for this class.) We consider several containment relations and study the cop number for classes
of graphs with a single forbidden graph with respect to these relations.

A graph is called \defn{$H$-minor-free} (\defn{$H$-topological minor-free}) if it does not contain $H$ as a minor (as a topological minor). 
Similarly, a graph is called \defn{$H$-subgraph-free} (\defn{$H$-free}) if it does not contain $H$ as a subgraph (as an induced subgraph). 

Families of graphs with unbounded cop number have been constructed \cite{AignerFromme84}. For every fixed $d\ge 3$, there even exist families of $d$-regular graphs
with unbounded cop number \cite{Andreae84}. On the other hand, Aigner and Fromme \cite{AignerFromme84} proved that the cop number of a planar graph is at most 3.
This result has been generalized to the class of graphs with genus $g$; Schroeder \cite{Schroeder} proved that the cop number of a graph is bounded by $\lfloor 3 g
/ 2\rfloor +3$.

Andreae~\cite{Andreae86} studied classes of $H$-minor-free graphs and showed that the cop number of a $K_5$-minor-free graph (or $K_{3,3}$-minor-free graph) is at most 3. Since a planar graph does not have a $K_5$ or $K_{3,3}$ as a minor this result extends the result on planar graphs. However, for our purposes the most interesting result of Andreae~\cite{Andreae86} is that for any graph $H$ the cop number is bounded in the class of $H$-minor-free graphs. In other words, forbidding one minor is enough to bound the cop number.

Andreae~\cite{Andreae86} also observed  
that excluding a topological minor does not necessarily bound
the cop number. In fact, it is an easy corollary from his work that the class of $H$-topological minor-free graphs has bounded cop number if and only if the
maximum degree of $H$ is at most 3.

Inspired by these results we study other containment relations:  
subgraphs and induced subgraphs. We give necessary and sufficient conditions for the class of
$H$-subgraph-free graphs ($H$-free graphs) to have bounded cop number. First we present our results for induced subgraphs.

\newcounter{thm:induced}
\setcounter{thm:induced}{\value{theorem}}
\begin{theorem}\label{thm:induced}
  The class of  $H$-free graphs has bounded cop number if and only if every connected component of $H$ is a path.
\end{theorem}

Let us remark that a single vertex is considered to be a path.
The graph consisting of a path on $\ell$ ($\ell \geq 1$) vertices is denoted by $P_\ell$. 
The backward implication of Theorem~\ref{thm:induced} is a consequence of the following
proposition.

\newcounter{prop:P-l-free}
\setcounter{prop:P-l-free}{\value{theorem}}
\begin{proposition}\label{prop:P-l-free}
  For every $\ell \ge 3$, every  $P_\ell$-free graph has cop number at most $\ell-2$.
\end{proposition}

Using the same technique, it is in fact possible to show the following stronger result:
\newcounter{prop:C-l-plus-free}
\setcounter{prop:C-l-plus-free}{\value{theorem}}
\begin{proposition}\label{prop:C-l-plus-free}
  For every $\ell \ge 3$, every  graph with no induced cycle of length at least $\ell$
   has cop number at most $\ell-2$.
\end{proposition}

Notice that it is possible to rephrase the condition of Theorem \ref{thm:induced} and say that every connected component of $H$ is a tree with at most two leaves.
Here is our result for $H$-subgraph-free graphs.

\newcounter{thm:subgraph}
\setcounter{thm:subgraph}{\value{theorem}}
\begin{theorem}\label{thm:subgraph}
  The class of  $H$-subgraph-free graphs has bounded cop number if and only if every connected component of $H$ is a tree with at most three leaves.
\end{theorem}

It is easy to see that the cop number of a tree is 1.
As an intermediate step towards Theorem~\ref{thm:subgraph},
we study how the cop number of a graph $G$ is related to its tree-width, which is denoted by 
$\tw(G)$.

\newcounter{thm:copno-treewidth}
\setcounter{thm:copno-treewidth}{\value{theorem}}
\begin{proposition}\label{thm:copno-treewidth}
  The cop number of a  graph $G$ is at most $\tw(G)/2 +1$.
\end{proposition}

This bound is sharp for tree-width up to $5$.

\section{Forbidding induced subgraphs}\label{sec:induced}

Our goal in this section is to prove Theorem~\ref{thm:induced}. Notice that a graph, whose every connected component is a path, is an induced subgraph of some
sufficiently long path. Hence, the following proposition proves the backward implication of Theorem~\ref{thm:induced}.

\setcounter{myremembertheorem}{\value{theorem}}
\setcounter{theorem}{\value{prop:P-l-free}}
\begin{proposition}
  For every $\ell \ge 3$, every  $P_\ell$-free graph has cop number at most $\ell-2$.
\end{proposition}
\setcounter{theorem}{\value{myremembertheorem}}
Let us remark that, for $\ell=1,2$, the cop number of a $P_{\ell}$-free graph is trivially 1.
\begin{proof}[Proof of Proposition~\ref{prop:P-l-free}]
  Let $G$ be a  $P_\ell$-free graph
  and let us also assume, without loss of generality, that $G$ is connected.  We will give a winning strategy for $\ell-2$ cops. Initially all $\ell-2$ cops are on
  the same arbitrary vertex.  The strategy is divided into stages. The distance between the cops and the robber is the minimum distance from the robber to a cop.
  The goal of each stage is to decrease the distance between the cops and the robber. Once the distance is decreased we begin the next stage. We will show that a
  stage lasts a finite number of rounds.

At the beginning of each stage we choose a \defn{lead cop} (for this stage) among the pieces which are at the minimum distance from the robber. All distances in this proof are measured after the robber's and before the cops' move. We route the lead cop and instruct the other pieces to follow the lead cop in single file; the cops should form a path of length $\ell-2$. 

If the distance between the cops and the robber is at most one, then the cops clearly win. Suppose that the distance between the lead cop on vertex $x$ and the robber on vertex $y$ is $d \ge 2$. We order the lead cop to travel along the shortest path from $x$ to $y$ and then follow the exact route the robber took from vertex $y$. Notice that since the graph is $P_\ell$-free  the distance between the cops and the robber will decrease after at most $\ell-d-1$ moves. Once the distance decreased, we move to the next stage.
\end{proof}

We mention the following result which can be derived using almost the same strategy as in Proposition~\ref{prop:P-l-free}.
\setcounter{myremembertheorem}{\value{theorem}}
\setcounter{theorem}{\value{prop:C-l-plus-free}}
\begin{proposition}
  For every $\ell \ge 3$, every  graph with no induced cycle of length at least $\ell$
   has cop number at most $\ell-2$.
    \qed
\end{proposition}

Before completing the proof of Theorem~\ref{thm:induced}, we look at bipartite graphs with no long induced paths. A simple modification of the proof of Proposition \ref{prop:P-l-free} yields a better bound for the bipartite case. Here is how the cops' strategy needs to be modified: the cops follow the lead cop in such a way that the distance between any two consecutive cops is 2. 
We leave the details of this proof to the reader.

\begin{proposition}\label{prop:P-l-free-bip}
  For every $\ell \ge 1$, every  $P_{2\ell}$-free bipartite graph has cop number at most $\ell$.
\end{proposition}


To prove the forward implication of Theorem~\ref{thm:induced}, we need to introduce two graph operations which do not decrease the cop number: clique substitution
and edge subdivision.  Let $N(v)$ be the the set of neighbors of a vertex $v$. A clique substitution at a vertex $v$ is to replace $v$ with a clique of size
$\lvert N(v)\rvert$ and create a matching between vertices of the clique and the vertices of $N(v)$. The graph obtained from a graph $G$ by substituting a clique
at each vertex of $G$ will be denoted by $G^+$. More formally, the vertex set of $G^+$ is $\bigcup_v ( \{v\}\times N(v) ) $ and two vertices $(v_1, u_1)$ and
$(v_2, u_2)$ are adjacent if and only if $v_1 = v_2$, or $v_1 = u_2$ and $u_1 = v_2$.

\begin{lemma}\label{lem:copno-clique}
  Clique substitution does not decrease the cop number.
\end{lemma}

\begin{proof}
  Let $\varphi\colon V(G^+)\to V(G)$ be the mapping such that $\varphi((v,u)) = v$ and let $k$ be the cop number of $G^+$.  Note that $\varphi^{-1}(v)$ is the
  vertex set of the clique which is substituted for $v$.  We simultaneously play two games: one on $G$ and another on $G^+$.  We assume that we have a winning
  strategy for the cop player on $G^+$ and we simulate her moves on $G$. On the other hand, the robber is playing on $G$ and we simulate his moves on~$G^+$.

  If $w^0_1,\dots,w^0_k$ are vertices on which cops are initially placed on $G^+$, then let $v^0_i:=\varphi(w^0_i)$ ($1\leq i \leq k$) be initial vertices of cops
  on $G$. Suppose the robber places his piece on vertex $r^0$ in $G$, we place the robber in $G^+$ on the vertex $s^0$ such that $\varphi(s^0)=r^0$.

  Let $r^0,r^2,r^4,\dots$ be the sequence of vertices visited by the robber on $G$. Clearly, for each $j$, either $r^j=r^{j+2}$ or $r^j$ is adjacent to $r^{j+2}$
  in $G$. Choose a sequence of vertices $s^0,s^1,s^2,\dots$ of $G^+$ with the following properties:

  \begin{enumerate}
    \renewcommand{\labelenumi}{(\alph{enumi})}
  \item $\varphi(s^{j})=\varphi(s^{j+1}) = r^{j}$ for all $j$,
  \item if $r^{j}=r^{j+2}$, then $s^{j}=s^{j+1}=s^{j+2}$,
  \item if $r^j$ is adjacent to $r^{j+2}$, then $s^{j+1}$ is the unique vertex in the clique $\varphi^{-1}(r^j)$ corresponding to $r_j$ which is adjacent to a
    vertex in $\varphi^{-1}(r^{j+2})$, and $s^{j+2}$ is the unique vertex in $\varphi^{-1}(r^{j+2})$ which is adjacent to $s^{j+1}$.
  \end{enumerate}
  
  \noindent It is easy to see that  $s^1,s^2,\dots$ is a walk in $G^+$.

  Now we let the $k$ cops play on $G^+$ against the robber whose sequence of moves is given by $s^1,s^2,\dots$.  The cops' positions are mapped to $G$ by
  $\varphi$.  Suppose that after a number of moves, the cops are on vertices $w^j_1,\dots,w^j_k$ of $G^+$.  Then, the robber moves from $r^j$ to $r^{j+2}$ on $G$.
  We let the cops on $G^+$ react first to the move $s^j$ to $s^{j+1}$ of the robber on $G^+$, this results in the cops' positions $w^{j+1}_1,\dots,w^{j+1}_k$.
  Then the move $s^{j+1}$ to $s^{j+2}$ results in cops' positions $w^{j+2}_1,\dots,w^{j+2}_k$.  Note that, by the construction of the graph $G^+$, for each cop
  $c$, the vertex $v^j_c:=\varphi(w^j_c)$ of $G$ is adjacent to the vertex $v^{j+1}_c:=\varphi(w^{j+2}_c)$.  Thus, it corresponds to a feasible cop move from
  $v^j_c$ to $v^{j+1}_c$ on $G$.

  Now suppose that the robber is captured in $G^+$.  This means that for a cop $c$, we have $s^j = w^j_c$.  If $j$ is even, this immediately implies that the cop
  is also captured in $G$.  However, if $j$ is odd, then there is a move $s^j\to s^{j+1}$ with $j$ even which we have to consider.  We let the cop $c$ follow the
  robber onto the vertex $s^{j+1}$ and thus catch the robber in $G$.
\end{proof}

The \defn{claw} is the complete bipartite graph with sides of size 1 and 3. The operation of clique substitution will be used to show that the cop number of claw-free graphs is unbounded.

\begin{lemma}\label{lem:copno-clawfree}
  The class of  claw-free graphs has unbounded cop number.
\end{lemma}

\begin{proof}
  Let $\mathcal{G}$ be a class of graphs with unbounded cop number and $\mathcal{G^+} := \{ G^+ \mid G \in \mathcal{G} \}$. Notice that all graphs in
  $\mathcal{G^+}$ are claw-free. Applying Lemma~\ref{lem:copno-clique}, we see that the cop number of graphs in $\mathcal G^+$ is unbounded.
\end{proof}

The other graph operation needed for the proof of Theorem \ref{thm:induced} is edge subdivision.  
Berarducci and Intrigila \cite{BerarducciIntrigila93} proved the following lemma.

\begin{lemma}[\cite{BerarducciIntrigila93}]\label{lem:copno-subdiv-even}
  Subdividing all edges of a graph an even number of times does not decrease
  the cop number.
\end{lemma}

This leads to the following result. Recall that the \defn{girth}
of a graph is the length of its shortest cycle if it has one, $+\infty$ otherwise.

\begin{lemma}\label{lem:copno-cyclefree}
  For every integer $\ell\ge 3$, the class of
   graphs with girth at least $\ell$ has unbounded cop number.
\end{lemma}
\begin{proof}

Let $\mathcal G$ be an arbitrary class of graphs with unbounded cop number.
For every $G \in \mathcal G$, let $G'$ be a graph with girth at least $\ell$ obtained from $G$ by subdividing all edges sufficiently often.  Let $\mathcal G' := \{ G' \mid G \in \mathcal G \}$. Applying Lemma~\ref{lem:copno-subdiv-even}, we see that the class $\mathcal G'$ has unbounded cop number.
\end{proof}

Now we are ready to complete the proof of Theorem~\ref{thm:induced}.

\setcounter{myremembertheorem}{\value{theorem}}
\setcounter{theorem}{\value{thm:induced}}
\begin{theorem}
    The class of  $H$-free graphs has bounded cop number if and only if every connected component of $H$ is a path.
\end{theorem}
\setcounter{theorem}{\value{myremembertheorem}}

\begin{proof}
The backward implication of the theorem follows from Proposition \ref{prop:P-l-free}. Indeed, notice that if every connected component of $H$ is a path, then $H$ is a subgraph of the path on $\abs H + p - 1$ vertices, where $p$ is the number of connected components of $H$. Hence, the cop number of an $H$-free graph is bounded by $\max\{|H| + p - 3, 1\}$.

Now we will prove the forward implication of the theorem. Let $H$ be a graph such that the class of $H$-free graphs has bounded cop number. Suppose that $H$ contains a cycle and let $\ell$ be the length of the longest cycle of $H$. Clearly, the class of graphs with no induced cycle of length at most $\ell$ is contained in the class of $H$-free graphs. However, by Lemma ~\ref{lem:copno-cyclefree} the class of graphs with no induced cycle of length at most $\ell$ has unbounded cop number; a contradiction. Hence, $H$ is a forest.

Now suppose that $H$ contains a vertex of degree at least $3$. Since $H$ is a forest, it must contain a claw as an induced subgraph. Clearly, the class of claw-free graphs is contained in the class of $H$-free graphs. However, by Lemma~\ref{lem:copno-clawfree} the class of claw-free graphs has unbounded cop number; a contradiction.  
Hence, $H$ is a forest of maximum degree at most 2, that is, $H$ is a disjoint union of paths.
\end{proof}

We note that in the second part of the proof (removing cycles) we could have used some known constructions which show that graphs simultaneously having an arbitrarily large girth and large cop number do exist; see for instance
Andreae~\cite{Andreae84} and Frankl~\cite{Frankl87dam}.

\subsection*{Some remarks about edge subdivisions.}
Lemma~\ref{lem:copno-subdiv-even} by Berarducci and Intrigila \cite{BerarducciIntrigila93} gives a bound on the cop number of graphs which result by uniformly
subdividing all edges an even number of times.  By modifying the proof of Lemma~\ref{lem:copno-clique}, the following can be shown.
\begin{lemma}\label{lem:copno-subdiv-once}
  Subdividing all edges of a graph does not decrease the cop number.
\end{lemma}

Combining Lemmas \ref{lem:copno-subdiv-even} and~\ref{lem:copno-subdiv-once} we obtain the general result.

\begin{corollary}\label{cor:copno-subdiv-general}
  For every positive integer $r$, subdividing every edge of a graph $r$ times does not decrease the cop number.
\end{corollary}

Berarducci and Intrigila \cite{BerarducciIntrigila93} noted that subdividing edges in a non-uniform manner can both increase and decrease the cop number.  As for
uniform subdivisions, however, it is possible to give an estimate.

\begin{proposition}\label{prop:subdiv-upperbound}
  Subdividing each edge $r$ times increases the cop number by at most one.
\end{proposition}
\begin{proof}[Proof (sketch)]
  Denote by $\widetilde G$ the graph which results from the graph $G$ by subdividing each edge $r$ times.  A winning strategy for $\cop(G)+1$ cops on $\widetilde
  G$ is the following.  Let an auxiliary cop pursue the strategy described for the lead cop in the proof of Proposition~\ref{prop:P-l-free}.  By this we make sure
  that the robber cannot change his direction or pass in the middle of a subdivided edge except for a finite number of times.  The other $\cop(G)$ cops simulate
  their winning strategy for $G$ on $\widetilde G$.
\end{proof}

To further enlighten what happens if edges are subdivided, we propose the following construction.  Let $G$ be an arbitrary graph with $n$ vertices and cop number
at least 2.  We construct a graph $\widehat G$ by adding paths of length $2n$ to $G$: every pair of non-adjacent vertices of $G$ is joined by such a path.  It is
not difficult to see that $\cop(\widehat G) = \cop(G)$.  But by subdividing edges of $\widehat G$, we can obtain a graph resulting from $K_n$ by subdividing every
edge $n$ times.  From Proposition~\ref{prop:subdiv-upperbound} we know that the cop number of this graph is at most~$2$.  

Considering this construction, it seems natural to propose the following conjecture, which implies the conjecture of Meyniel (see Frankl~\cite{Frankl87dam}) that
$\cop(G)$ is in $O(\sqrt{\abs{G}})$.

\begin{conjecture}
  For graphs $G$ obtained by subdividing edges of complete graphs $K_n$ we have $\cop(G)$ in $O(\sqrt n)$.
\end{conjecture}

\section{Forbidding (not necessarily induced) subgraphs}\label{sec:subgraphs}

We now turn our attention to classes of graphs for which we forbid (not necessarily induced) subgraphs.
One key ingredient for the proof of Theorem~\ref{thm:subgraph} is the fact that families of graphs with bounded circumference have bounded cop number.  Although
this already follows from Proposition~\ref{prop:C-l-plus-free}, in this section we give a better upper bound based on an estimate on the cop number in terms of the
tree-width, which we believe to be of interest in its own.

Let us first briefly recall the definition of the tree-width of a graph.  A \defn{tree decomposition} of a graph $G$ is a pair $(T, \{W_{x} \mid x \in V(T)\})$
where $T$ is a tree, and $\{W_{x} \mid x \in V(T)\}$ a family of subsets of $V(G)$ (called ``bags'') such that
\begin{itemize}
\item $\bigcup_{x \in V(T)} W_{x} = V(G)$;
\item for every edge $uv \in E(G)$, there exists $x \in V(T)$ with $u,v  \in W_{x}$, and
\item for every vertex $u\in V(G)$, the set $\{x\in V(T) \mid u \in W_{x}\}$ induces a subtree of $T$.
\end{itemize}
The \defn{width} of tree decomposition $(T,  \{W_{x} \mid x \in V(T)\})$ is
$\max \{ |W_{x}| - 1 \mid x \in V(T)\}$.
The \defn{tree-width} $\tw(G)$ of $G$ is the minimum width among all tree decompositions
of $G$. We refer the reader to Diestel's book~\cite{Diestel-GTIII} for an introduction
to the theory around tree-width.

Our proof of Proposition~\ref{thm:copno-treewidth} relies on a well-known
strategy for the cops and robber game: guarding a shortest path. Assume that $P$ is a shortest $uv$-path,
for two distinct vertices $u,v$ of a graph $G$, and that a cop is sitting at the beginning
on some vertex of $P$. The cop's strategy consists in moving along $P$ in such a way
that his distance to $u$ is as close as possible to the robber's distance to $u$.
It is easily seen that, after a finite number of initial moves, when it is the robber's turn to play,
the cop's distance to $u$ will be the same as the robber's distance to $u$ when the latter is no
more than $|P|$. This ensures that the robber cannot go on any vertex of $P$ without
being caught. (This strategy has been first used by Aigner and Fromme \cite{AignerFromme84},
in their proof that the cop number of planar graphs is at most 3.)

\setcounter{myremembertheorem}{\value{theorem}}
\setcounter{theorem}{\value{thm:copno-treewidth}}
\begin{proposition}
  The cop number of a  graph $G$ is at most $\tw(G)/2 +1$.
\end{proposition}
\setcounter{theorem}{\value{myremembertheorem}}

\begin{proof}
  Let us consider an optimal tree decomposition of $G$. Since the
  tree-width of $G$ equals the maximum tree-width of its connected components, we may assume
  without loss of generality that $G$ is connected.
  For a bag $X \subseteq V(G)$ of the tree decomposition, we denote by
$t_{X}$ the vertex of $T$ corresponding to $X$.

At the beginning, an arbitrary bag $B \subseteq V(G)$ of the tree decomposition is selected, and
all its vertices are guarded in the following way: Letting $b_{1}, b_{2}, \dots, b_{k}$
denote the vertices in $B$, we let the $i$th cop ($1 \leq i \leq \lfloor k  / 2 \rfloor$)
guard a shortest $b_{2i - 1}b_{2i}$-path in $G$, and, if $k$ is odd, we put an additional
cop on vertex $b_{k}$. This ensures that, after a finite number of moves, the robber cannot
go on any vertex in $B$, and hence is confined to (the subgraph corresponding to)
some tree $T'$ of $T \setminus t_{B}$.
(We may assume $B \neq V(G)$, as otherwise the robber is trivially caught.)

Let $B' \subseteq V(G)$ be the unique bag of the tree decomposition that is adjacent to $B$ in $T$
with $B' \cap C \neq \varnothing$. 
Observe that $B \cap B'$ is a cutset of the graph $G$.
We show that the cops can move in such a way
that the vertices of $B \cap B'$ remains guarded, and after a finite number of moves
all the vertices of $B'$ (instead of $B$) are guarded.

Consider each cop.
Suppose first that the cop sits on a vertex of $B  \setminus B'$
or guards a shortest path between two vertices in $B  \setminus B'$. Then
we send him to guard a shortest path between two unguarded vertices in
$B' \setminus B$ (or to sit on the last unguarded vertex if there is only one such vertex).
Assume now that the cop sits on a vertex of $B \cap B'$ or
guards  a shortest $b_{i}b_{j}$-path with $b_{i} \in B \cap B'$ and $b_{j} \in B \setminus B'$.
Then the cop first goes to $b_{i}$ (if he is not already there) along the path he keeps.
The he starts guarding an
arbitrary $b_{i}b'_{j'}$-path, where $b'_{j'}$ is any unguarded vertex of $B' \setminus B$.
Notice that, while it may take some moves before all the vertices of the path are safely guarded,
at least the vertex $b_{i}$ is guarded at every time.
Suppose finally that the cop guards a shortest $b_{i}b_{j}$-path with $b_{i}, b_{j} \in B \cap B'$.
In this case, the cop does not modify his strategy, and keeps guarding his path.

After a finite number of moves all the vertices in $B'$ are guarded,
and the robber did not have, at any time, the opportunity
to go on a vertex in $B \cap B'$ without being caught.
Moreover, the number of necessary cops is at most $\lceil |B'| / 2 \rceil \leq \tw(G)/2 + 1$,
and the robber is reduced to stay in
(the subgraph corresponding to) some tree of $T \setminus t_{B'}$ which is a proper
subtree of $T'$. Therefore,
by repeating this operation a finite number of times the robber will eventually be caught.
This completes the proof.
\end{proof}

We remark that the bound given in
Proposition~\ref{thm:copno-treewidth} is best possible for small values of the
tree-width: For every $k=1, 2, \dots, 5$, there are graphs with tree-width $k$
and cop number $\lfloor k/2 \rfloor + 1$ (this is easily seen for $k=1,2,3$, and
the Petersen graph and the graph which is the disjoint union 
of the Petersen graph and a complete graph on 6 vertices
are such examples for $k=4$ and $5$, respectively).
On the other hand, we do not know whether there exists a constant $c > 0$ and
an infinite family of graphs such that $\cop(G) \geq c\cdot \tw(G)$ holds for every graph $G$
in the family.

Let us recall that the \defn{circumference} of a graph is the length of its
longest cycle if it has one, $+\infty$ otherwise.

\begin{corollary}\label{cor:copno-circumference}
  The cop number of a  graph is less than
  or equal to half its circumference. 
\end{corollary}
\begin{proof}
  It is a well-known fact that $\tw(G) \le \circum(G)-1$
  holds for every graph $G$, where $\circum(G)$ denotes the circumference of  $G$ 
  (see for instance Exercise~12.18 in Diestel's book~\cite{Diestel-GTIII}).  With Proposition~\ref{thm:copno-treewidth}, we conclude $\cop(G) \le
  \circum(G)/2$.
\end{proof}

\setcounter{myremembertheorem}{\value{theorem}}
\setcounter{theorem}{\value{thm:subgraph}}
\begin{theorem}
  The class of  $H$-subgraph-free graphs has bounded cop number if and only if every connected component of $H$ is a tree with at most three leaves.
\end{theorem}
\setcounter{theorem}{\value{myremembertheorem}}
\begin{proof}
We first show that the requirements in the statement of the theorem are necessary. Let $H$ be a graph such that the family $\mathscr F$ of connected
graphs not containing $H$ has bounded cop number.

  First, suppose that $H$ contains a cycle, and let $\ell$ be the length of a longest cycle in $H$.  Then
  $\mathscr F$ contains the family of connected graphs with girth at least $\ell +1$.
  However, by Lemma~\ref{lem:copno-cyclefree}, the cop number of this family is unbounded.
  Hence, $H$ is a forest.

  Second, suppose that $H$ has a vertex of degree at least 4.
  This implies that $\mathscr F$ contains all connected graphs with maximum degree 3, but
  Andreae~\cite{Andreae84} proved that there exists a family of 3-regular graphs on which the cop number is unbounded.  Hence, $H$ has maximum degree at most 3.

  Third, suppose that there is a tree in $H$ which has two vertices of degree 3.
  Let $\ell$ denote the distance between these two vertices in $H$.  Now
  $\mathscr F$ contains the family of all those connected graphs in which
  every two vertices of degree 3 or more have distance at least $\ell+1$.  Starting
  from an arbitrary family of graphs on which the cop number is unbounded, a family with this property can be constructed by subdividing every edge $\ell$ times,
 as follows from Corollary~\ref{cor:copno-subdiv-general}.
 Thus, each connected component of $H$ contains at most one vertex of degree 3.

  We now show that any $H$ meeting the conditions in the theorem yields a family of graphs with bounded cop number.
  The proof will be by induction on the number of connected components of $H$.  For a single component, by Proposition~\ref{prop:P-l-free}, we may assume that a vertex of
  degree 3 does in fact exist.
  We will prove the following claim.

\medskip

  \noindent \textbf{Claim.}  Let $H$ be a tree with maximum degree 3 which has precisely one vertex $v$ of degree 3.  Denote by $r$ the maximum distance of a
  vertex from $v$.  If $G$ does not contain $H$,  then $\cop(G) \le 2r$.

\medskip

  Before we prove the claim, let us complete the induction.  The start of the induction is settled.  Let $T$ be a connected component of $H$, and assume that
  $\cop(G') \le k$ for every graph $G'$ not containing $H\setminus V(T)$.  
  Let $G$ be a graph not containing $H$.  If $G$ does not contain $T$, we are done by the
  claim and the remark preceding it.  Otherwise, let $T'$ be a subgraph of $G$ isomorphic to $T$.  We place $\abs T$ cops on the vertices of $T'$.  This corners
  the robber in a connected component of $G\setminus V(T')$.  
  Noting that $G\setminus V(T')$ does not contain $H\setminus V(T)$, by induction, by restricting to the
  connected component containing the robber, 
  we can catch the robber in $G\setminus V(T')$ using $k$ cops.  
  This bounds the cop number of $G$ by $k+\abs T$ and
  concludes the induction.

  \textbf{Proof of Claim.}  We prove the claim in the case when each leaf of $H$ has distance exactly $r$ from $v$.  The general case follows easily from this.

  By Proposition~\ref{prop:P-l-free}, we may assume that $G$ contains a path $P$ on $2r$
  vertices, because otherwise we have $\cop(G) \le 2r$. 
  We guard the path by placing $r$ cops on every other vertex of $P$, 
  and show that what remains of $G$ has cop number at most $r$.  Assume that $G\setminus V(P)$
  contains a cycle $C$ of length at least $2r+1$.  
  Then, since $G$ is connected, we can identify a subgraph isomorphic to $H$ choosing $v$ to be a vertex on $C$
  which has minimum distance to a vertex in $P$, while two of the three branches of the tree are wound around $C$, the other extends to $P$.  Hence, $G\setminus
  V(P)$ contains no such cycle.  By invoking Corollary~\ref{cor:copno-circumference} for the connected component of $G\setminus V(P)$ containing the robber, we see that the
  cop number of $G\setminus V(P)$ is at most $r$.
\end{proof}

%

\input{bibliography}
\end{document}

%% file: defs.amsart.tex
\theoremstyle{plain}
\newtheorem{theorem}{Theorem}
\newtheorem{proposition}[theorem]{Proposition}
\newtheorem*{proposition*}{Proposition}
\newtheorem{corollary}[theorem]{Corollary}
\newtheorem{lemma}[theorem]{Lemma}

\newtheorem*{theorem*}{Theorem}
\newtheorem*{lemma*}{Lemma}
\newtheorem*{conjecture*}{Conjecture}
\newtheorem{conjecture}[theorem]{Conjecture}

\theoremstyle{definition}

\newtheorem*{exercise*}{Exercise}

\theoremstyle{remark}

\newtheorem*{remark*}{Remark}
\newtheorem{remsTh}[theorem]{Remarks}

\newcommand{\subclass}[1]{}

\newcommand{\enumTi}[1]{\renewcommand{\theenumi}{#1}}

\input{defs.common}


%% file: defs.common.tex
\renewcommand{\em}{\sl}

\newcommand{\lt}{\left}
\newcommand{\rt}{\right}

\newcommand{\abs}[1]{{\lt\lvert{#1}\rt\rvert}}

\newlength{\algotabbingwidth}
